\documentstyle[12pt,amssymb]{amsart}

\def\eop{\blacksquare}

   \def\B{\cal B} \def\S{\cal S}
 \def\G{\Gamma}  \def\Y{{\operatorname Y}}

\def\k{\operatorname{\bf k}}

 \def\X{\bf X} 
\def\Nat{\operatorname{Nat}}

 \def\group{\operatorname{group}}
 \def\Bn{\operatorname{B}_{n}}
 
 \def\Hom{\operatorname{Hom}}
 \def\End{\operatorname{End}}
 
\def\I{\operatorname{I}}  
  \def\Gt{{\G}^{\otimes}}
\def\TG{\operatorname{T\G}}  
 \def\m{\operatorname{m}} \def\tr{\triangle}
 
\def\SLA{\operatorname{SLA}} \def\SLC{\operatorname{SLC}}
\def\LA{\operatorname{LA}} \def\LC{\operatorname{LC}}

\newcommand{\im}{\operatorname{im}} 
\newcommand{\id}{\operatorname{id}} \newcommand{\lin}{\text{lin}}
 
 \newcommand{\eid}{\operatorname{id}}

 \newcommand{\der}{\text{der}}
\newcommand{\ben}{\begin{enumerate}} \newcommand{\de}{\operatorname{d}}
\newcommand{\deh}{\hat{d}}

 \newcommand{\een}{\end{enumerate}}

\newtheorem{LEM}[DEF]{Lemma} \newtheorem{TH}[DEF]{Theorem}
 \theoremstyle{definition}
\newtheorem{DEF}{Definition}[section]

\begin{document} \title{Koszul theorem for S-Lie coalgebras} \author{Jerzy
R\'o\.za\'nski} \thanks {Supported by KBN grant \# 2 P302 023 07}
\date{December 1995} \address{Department of Mathematical Methods in Physics,
University of Warsaw, Ho\.za 74, PL 00-682 Warszawa, Poland}
\email{rozanski@@fuw.edu.pl} \maketitle

\begin{abstract} For a symmetry braid $\S$-Lie coalgebras, as a dual object
to algebras introduced by Gurevich \cite{DG}, are considered.  For an Young
antisymmetrizer $Y(-\S)$ an $\S$-exterior algebra is introduced.  From this
differential point of view $\S$-Lie coalgebras are investigated.  The dual
Koszul theorem in this case is proved.  \end{abstract}

\section{Introduction} For a symmetric braid $\S$ we wish to investigate
$\S$-Lie coalgebras as a dual object to $\S$-Lie algebras introduced by
Gurevich \cite{DG}.  With an $\S$-Lie coalgebra we consider an $\S$-exterior
algebra over a coalgebra given by an Young antisymmetrizer  $\Y (-\S)$
\cite{Boz}, \cite{MR}.  The aim of this paper is to find out a kind of the
Koszul theorem \cite{Kosz} for $\S$-Lie coalgebras. This means that the
Jacobi condition for an $\S$-Lie coalgebra is equivalent to the differential
condition $\deh^{2}=0$ for a derivation of an $\S$-exterior algebra. We are
interested in finding the cohomological meaning of the Jacobi condition for
$\S$-Lie (co)algebras \cite{DG}.

In section 1 we review definitions in a braided monoidal category $\X$.  In
the paper we use the following notations.  Let $\k$ be a field and $\G$ be a
$\k$-space, a dual $\k$-space $\G^{\ast}=\Hom (\G,\k)$ and $\TG$ means a
tensor algebra for $\k$-space $\G$.  We consider a strict monoidal category
$\X$ generated by one object $\G\in\ Obj (\X)$, and $Mor (\X)$ is a set of
morphisms.  By $\de\in\der (\TG)$ let us denote a set of derivations of a
free graded algebra $\TG$.

For the symmetry braid $\S$ we have a symmetric monoidal category and the
definition of a Lie (co)algebra is given over this category. For the braid
group we consider the notion of an Young antisymmetrizer.  As an example we
have an $\S$-antisymmetrizer $\Y (-\S)$.  An $\S$-exterior algebra is a
factor algebra of a free algebra $\TG$ and an ideal is equal to a kernel of
the operator $\Y (-\S)$.

Section 2 contains a brief summary of definitions and remarks of a Lie
(co)algebra over a monoidal category. We will touch only a few aspects of the
theory of Lie coalgebras. for a more complete theory see \cite{WM}.  New is
the dual Koszul theorem for Lie coalgebras. This is possible, because we show
that for any derivation of a free tensor algebra exists a factor derivation
of an exterior algebra.

In section 3 we consider the notion of $\S$-Lie coalgebras over an
$\S$-braided monoidal category. The basic idea is that the morphism condition
for a comultiplication guarantes the existence of a factor derivation of a
free tensor algebra factor by the kernel of the Young antisymmetrizer
$\Y(-\S)$.

The main result of this paper is the dual Koszul theorem for $\S$-Lie
coalgebras as the possible form of the Jacobi identity \cite{OPR}.  Let
$\bigwedge_{\S} \G=\Gt/\ker Y(-\S)$ and $\deh\in\der (\TG/\ker \Y(-\S))$,
$\S_{1}=\S\otimes\id$ and $\S_{2}=\id\otimes\S$.  For an $\S$-exterior
algebra ${\cal E}_{\S} (\G)=(\bigwedge_{\S} \G,\;\deh_{\tr})$ over a $\S$-Lie
coalgebra $\SLC=(\G,\;\tr,\;\S)$ holds the equation $\deh_{\tr}|\G=\tr$.  \\
\\ \noindent {\bf Theorem.}{(\bf Koszul Theorem for $\S$-Lie coalgebras}).\\
For the Young antisymmetrizer $\Y (-\S)$ let ${\cal E}_{\S} (\G)$ be the
$\S$-exterior algebra over the  $\S$-Lie coalgebra $\SLC$.  Then the
following conditions are equivalent \begin{itemize} \item
$\deh^{2}=0\in\G^{\wedge 3},$ \item $\im (\de_{\tr})^{2}\subset\ker Y_{3}
(-\S)$, \item
$(\id+\S_{1}\S_{2}+\S_{2}\S_{1})\circ(\tr\otimes\id)\circ\tr=0$, \item
$(\tr\otimes\id)\circ\tr - (\id\otimes\tr)\circ\tr -
(\id\otimes\S)\circ(\tr\otimes\id)\circ\tr=0.$ \end{itemize}

\section{Notations in Braided monoidal category} In this section we recall
some definitions and remarks for a braided monoidal category \cite{SM},
\cite{DY}.  A category $\X$  equipped with an object $\I$ and a bifunctor
$\otimes:\X\times\X\rightarrow\X$, which is associative and for which $\I$ is
a twosided identity, we call {\em a monoidal category}.  In this paper
objects are linear spaces, an identity $\I$ is a field $\k$ and bifunctors
$\otimes$ means $\otimes_{\k}$.  We consider only a strict monoidal category
$\X$ and denote a set of morphisms for this category by $Mor (\X)$.  A
natural transformation $\B\in\Nat(\otimes,\;\otimes^{opp})$ is a family of
morphisms \begin{equation} \B_{X,Y}\in\lin(X\otimes Y,Y\otimes X),
\end{equation} which for all morphisms $f\in\Hom(X,Y)$ and $g\in\Hom(Z,W)$
satisfies the naturality conditions:  \begin{equation} (g\otimes
f)\circ\B_{X,Z}=\B_{Y,W}\circ(f\otimes g).  \end{equation}

        In Particular $\B_{V,W}$ is $\B$-morphism if we have two tetragons
\begin{gather} (\id\otimes\B_{V,W})\circ\B_{V\otimes
W,U}=(\B_{V,W}\otimes\id)\circ\B_{W\otimes V,U},\nonumber\\
(\B_{V,W}\otimes\id)\circ\B_{V\otimes
W,U}=(\id\otimes\B_{V,W})\circ\B_{W\otimes V,U}, \end{gather}

{\bf Prebraid.} A natural transformation $\B\in\Nat(\otimes,\;\otimes^{opp})$
is a {\em prebraid} if two trigons hold:  \begin{eqnarray} {\cal B}_{V\otimes
W,U}&=&({\cal B}_{V,U}\otimes\eid_W) \circ(\eid_V\otimes{\cal
B}_{W,U}),\nonumber\\ {\cal B}_{V,W\otimes U}&=&(\eid_W\otimes{\cal B}_{V,U})
\circ({\cal B}_{V,W}\otimes\eid_U), \end{eqnarray} The consequence of the
	naturality condition  and pairs of trigons is the braided hexagon

{\bf Prebraided monoidal category.} A monoidal category $\X$ is {\em
prebraided} if it is equipped with the prebraid $\B\in
Nat(\otimes,\otimes^{opp})$.

        A prebraided monoidal category is {\em braided} if $\B$ is a natural
isomorphism. Then we have pairs of braidings:  \begin{equation}
\B\in\Nat(\otimes,\otimes^{opp})\qquad
{\text{and}}\qquad\B^{-1}\in\Nat(\otimes^{opp},\otimes), \end{equation} such
that $\B\circ\B^{-1}=\id_{\otimes^{opp}},\;\B^{-1}\circ\B=\id_{\otimes}$,
$(B_{U,V})^{-1}=(\B^{-1})_{V,U}$.

        A braided monoidal category is {\em symmetric} if the braid $\S$
satisfies the condition:  \begin{equation}
\S_{U,\;V}\circ\S_{V,\;U}=\id_{V\otimes U}.  \end{equation} This means that
the notion $\S=\S^{-1}$ is badly posed  from the categorical point of view.

For the symmetry braid equalled to the flip $\S=\tau$ \begin{gather} \tau_{U,
V }(u\otimes v)=v\otimes u,\qquad u\in U,\; v\in V.  \end{gather}
\begin{LEM}\label{taumor} Any morphism is a $\tau$-morphism in a
$\tau$-monoidal category.  \end{LEM}

{\bf Young antisymmetrizer.} Let $\Bn$ be a braid group with generators\\
$\{\sigma_{1},\sigma_{2},\ldots,\sigma_{n-1}\}$ and relations \begin{gather}
\sigma_{i}\sigma_{j}=\sigma_{j}\sigma_{i}\qquad |i-j| > 1,\nonumber\\
\sigma_{i}\sigma_{i+1}\sigma_{i}=\sigma_{i+1}\sigma_{i}\sigma_{i+1}\qquad
i=1,\;2,\ldots,\;n-2.  \end{gather}

Let $\B$ be the Yang Baxter operator, i.e. an invertible endomorphism of the
$\G^{\otimes 2}$ which satisfies the braid equation \begin{equation}
(\id\otimes \B)(\B\otimes \id)(\id\otimes \B)= (\B\otimes\id)(\id\otimes
\B)(\B\otimes\id)\in\End(\G^{\otimes 3}).  \end{equation} Let
$B_{k}=\id_{k-1}\otimes B\otimes\id_{n-k-1}\in\End (\G^{\otimes n})$ be the
set of endomorphisms of $\G^{\otimes n}$.  For the Yang Baxter operator $\B$
we can define the representation $\rho_{\B}$ of the braid group $B_{n}$ in
the $\k$-space $\G^{\otimes n}$ \begin{equation}
\rho_{\B}\in\group\{B_{n},\;\End(\G^{\otimes n})\}:\qquad
\rho_{\B}(\sigma_{k})=B_{k}.  \end{equation}

Let $\psi$ mean the injection map from a permutation group $P_{n}$ to the
braid group $B_{n}$ \cite{SW}. Then we have the image of the map $\psi$, the
subset $\Xi_{n}$ \begin{equation} \Xi_{n}=\psi (P_{n})\subset B_{n}.
\end{equation}

For an Yang Baxter operator $\B$, the Young (braided) antisymmetrizer
$\Y(\B)$ is defined by \cite{Boz}, \cite{MR} \begin{equation}
\Y(\B)=\sum_{b\in\Xi_{n}} \rho_{\B}(b).  \end{equation}

For $\B=-\S$ we have the Woronowicz form of the braided antisymmetrizer with
a sign of a permutation \cite{SW} \begin{equation} \Y(-\S)=\sum_{b\in\Xi_{n}}
sign[\psi^{-1} (b)]\; \rho_{\S}(b).  \end{equation}

\section{Lie coalgebras over monoidal category} Michaelis \cite{WM} defined a
Lie algebra and a Lie coalgebra over a monoidal category. We recall these
definitions and the Koszul theorem in the dual form is rewritten.
\begin{DEF}[Lie algebra, \cite{WM}] A Lie algebra $\LA$ over a monoidal
category $\X$ is a pair $\{\G,\;\m\}$, where $\G\in Obj(\X)$ and  a map
$\m\in Mor(\X)$ is subjected to \begin{itemize} \item $\tau\circ
(\m\otimes\id) =(\id\otimes \m)\circ \tau_{1}\tau_{2}$, \item
$J(\m)\equiv\m\circ (\m\otimes\id)\circ
(\id+\tau_{1}\tau_{2}+\tau_{2}\tau_{1})=0$, \item
$\ker\Y_{2}(-\tau)\subset\ker\m$.  \end{itemize} \end{DEF} Remarks.  The
first condition is the $\tau$-morphism condition for the multiplication.
Usually it is omitted, because it is trivial $\forall \m$, see lemma
\ref{taumor}. We write this for the pedagogical reasons.  The second is the
Jacobi identity for the multiplication in the argument free form.  For the
map $J(\m)$ we have $J(\m):\G^{\wedge 3}\rightarrow\G$. The third condition
is the commutativity of the multiplication $\m\circ (\id+\tau)=0$.

For the Jacobi condition in the free argument form $J(\m)=0$ we have another
possible forms.  \begin{LEM} The following assertions are equivalent
\begin{itemize} \item $J(\m)\equiv\m\circ (\m\otimes\id)\circ
(\id+\tau_{1}\tau_{2}+\tau_{2}\tau_{1})=0$, \item $m\circ
(\m\otimes\id)-\m\circ (\id\otimes\m)-m\circ (\id\otimes\m)\circ
(\id\otimes\tau)=0,$ \item $\forall
x,\;y,\;z\in\G:\;\;[[x,\;y],\;z]-[x,\;[y,\;z]]-[[x,\;z],\;y]=0.$
\end{itemize} \end{LEM}

\begin{DEF}[Lie coalgebra, \cite{WM}] A Lie coalgebra $\LC$ over a monoidal
category $\X$ is a pair $\{\G,\,\tr\}$, where $\G\in Obj(\X)$ and $\tr\in
Mor(\X)$ is subjected to \begin{itemize} \item
$(\tr\otimes\id)\circ\tau=\tau_{2}\tau_{1}\circ(\id\otimes\tr)$, \item
$cJ(\m)\equiv
(\id+\tau_{1}\tau_{2}+\tau_{2}\tau_{1})\circ(\tr\otimes\id)\circ\tr=0,$ \item
$\im \tr\subset\im Y_{2}(-\tau)$.  \end{itemize} \end{DEF} Remarks.  The
first condition is the $\tau$-morphism condition for the comutliplication. It
is omitted, because it is trivial $\forall \tr$ see lemma \ref{taumor}.  The
second is the Jacobi identity for the comultiplication $\tr$ in the argument
free form.  For the co-Jacobiator $cJ(\m)$ we have
$cJ(\m):\G\rightarrow\G^{\wedge 3}$.  The third condition is the
cocommutativity of the comultiplication $(\id+\tau)\circ\tr=0$.

Consider an exterior algebra ${\cal E}_{\tau}
(\G^{\ast})=\{\bigwedge\G^{\ast},\;\deh_{\tr}\}$ over a Lie algebra
$\LA=\{\G,\;\m\}$. This means that the equation
$\deh_{\tr}|\G^{\ast}=\m^{\ast}$ holds.  \begin{TH}[Koszul Theorem,
\cite{Kosz}]\label{Kthm} For an exterior algebra ${\cal E}_{\tau}
(\G^{\ast})$ over Lie algebra $\LA$ the following assertions are equivalent
\begin{itemize} \item $J(\m)=0\in\G,$ \item $(\deh_{\tr})^{2}=0\in\G^{\ast
\wedge 3}.$ \end{itemize} \end{TH} {\em Proof.} Let $\alpha\in\G^{\ast}$ and
$x,\; y,\; z\in\G$.  Let us introduce two operators $e_{x}y=x\wedge
y,\;\;\;(i_{x}\alpha)y=\alpha (e_{x} y).$ Then the Lie derivation is ${\cal
L}_{X}=\deh\circ i_{X}+i_{X}\circ\deh.$ For $(\deh_{\tr})^{2}=0$ we have
\begin{multline*} 0=<(\deh)^{2} \alpha,\; x\wedge y\wedge z>=<i_{X}\deh^{2}
\alpha,\; y\wedge z>\\ =<(i_{X}\circ\deh)\deh \alpha,\; y\wedge z>= <({\cal
L}_{X}-\deh\circ i_{X})\deh \alpha,\; y\wedge z>\\ =-<\deh\circ i_{X})\deh
\alpha,\; y\wedge z> +<\alpha,\; {\cal L}_{X} (y\wedge z)>.  \end{multline*}
Consider two equations \begin{multline*} <\deh\circ i_{X})\deh \alpha,\;
y\wedge z>=<i_{X})\deh \alpha,\; \m(y\wedge z)>=<\deh \alpha,\; e_{X}
\m(y\wedge z)>\\ =<\alpha,\; \m(x\wedge \m(y\wedge z)>=<\alpha,\;
[x,[y,z]],\\ <\alpha,\; {\cal L}_{X} (y\wedge z)>=<\alpha,\; ({\cal L}_{X}
(y))\wedge z+y\wedge ({\cal L}_{X} (z))>\\ =<\deh\alpha,\;[x,y]\wedge
z+y\wedge [x,z]>=<\alpha,\;[[x,y],z]+[y,[x,z]]>.  \end{multline*} Then
$$<\deh^{2} \alpha,\; x\wedge y\wedge
z>=<\alpha,\;[[x,y],z]+[y,[x,z]]-[x,[y,z]]>=0.\eop$$

For a derivation of the free tensor algebra $\TG$ we have a factor derivation
$\deh\in \TG/\ker \Y(-\tau)$.  \begin{LEM}[Factor derivation]\label{fdlie}
The factor derivation $\deh_{\tr}\in\der(\TG/\ker\Y(-\tau))$ exists for any
derivation $\de_{\tr}$ of a free algebra $\TG$.  \end{LEM} {\em Proof.} This
follows from the fact that any comultiplication $\tr$ is the morphism in the
$\tau$-monoidal category.  Using the recurrent formula for the derivation
$\de_{\tr}$ $$\de_{\tr}|\G^{\otimes n}=\sum_{k=1}^{n-1}
({-1}^{k})\id_{k-1}\otimes\tr\otimes\id_{n-k-1}.$$ we can proof this fact by
induction. $\eop$

We can consider an exterior algebra ${\cal E}_{\tau}
(\G)=(\bigwedge\G,\;\deh_{\tr})$ over a Lie coalgebra $\LC=(\G,\;\tr)$. Then
the equation $\deh_{\tr}|\G=\tr$ holds.  \begin{TH}[Koszul Theorem for Lie
coalgebras]\label{dKthm} For an exterior algebra ${\cal E}_{\tau} (\G)$ over
a Lie coalgebra $\LC$ the following assertions are equivalent:
\begin{itemize} \item $\deh_{\tr}^{2}=0\in\G^{\wedge 3}$, \item
$\im(\de_{\tr})^{2}\subset\ker\Y(-\tau),$ \item
$cJ(\tr)=(\id+\tau_{1}\tau_{2}
+\tau_{2}\tau_{1})\circ(\tr\otimes\id)\circ\tr=0,$ \item
$(\tr\otimes\id)\circ\tr-(\id\otimes\tr)\circ\tr
-(\id\otimes\tau)\circ(\tr\otimes\id)\circ\tr=0.$ \end{itemize} \end{TH} {\em
Proof.} The first and the second conditions are equivalent from the lemma
\ref{fdlie} $$\Y_{3} (-\tau)\circ \de_{\tr}^{2}|\G=0\Leftrightarrow
\deh_{\tr}^{2}=0.$$ The remark is that for a derivation we can check
conditions on generators of the algebra. From the identity in lemma
\ref{blumen} the second condition is equivalent to the fourth condition. From
this condition by the commutativity of the multiplication $\m\circ\tau=-\m$
we get the third condition.  $\eop$

Remark.  The third condition is the Jacobi formula for the comultiplication
$\tr$ in the Woronowicz form \cite{SW}.

\section{S-Lie coalgebras} Let us recall the definition of the S-Lie algebras
introduced by D. Gurevich \cite{DG}. In this section $\S\in\End (\G^{\otimes
2})$ is a symmetry braid, $\S^{2}=\id$. A category $\X$ is $\S$-braided.
\begin{DEF}[$\S$-Lie algebra, \cite{DG}] An S-Lie algebra $\SLA$ over
$\S$-braided category $\X$ is a pair $(\G,\;m)$, where $\G\in Obj(\X)$ and
the multiplication $\m\in Mor (\X)$ is subjected to \begin{itemize} \item
$\S\circ(\id\otimes\m)=(\m\otimes\id)\circ\S_{1}\S_{2}$ \item
$\m\circ(\m\otimes \id)\circ(\id+\S_{1}\S_{2}+\S_{2}\S_{1})=0$, \item
$\ker\m\subset\ker Y_{2} (-\S)$ \end{itemize}\end{DEF} Remarks.  For the
first condition if $\S$ is the symmetry braid then the following morphism
conditions for the multiplication are equivalent $$ \S\circ (\m\otimes
\id)=(\id\otimes \m)\circ \S_{2}\S_{1}\Leftrightarrow \S\circ (\id\otimes
\m)=(\m\otimes \id)\circ \S_{1}\S_{2}.  $$ The first condition is the strong
condition for multiplications.  For example \begin{equation}\label{color}
\S(e_{i}\otimes e_{j})=\epsilon_{ij}\cdot e_{j}\otimes e_{i},\qquad
\epsilon_{ij}\epsilon_{ji}=1,\qquad \S^{2}=\id.  \end{equation} Then the
comultiplication $\tr$ is $\S$-morphism for color algebras \cite{Sch}.

        \begin{DEF}[$\S$-Lie coalgebra] An S-Lie coalgebra $\SLC$ over
$\S$-braided monoidal category $\X$ is a pair $(\G,\; \tr)$, where $\G\in Obj
\X$ and the comultiplication $\tr\in Mor (\X)$ is subjected to
\begin{itemize} \item
$(\tr\otimes\id)\circ\S=\S_{1}\S_{2}\circ(\id\otimes\tr).$ \item
$cJ(\tr)=(\id+\S_{2}\S_{1}+\S_{1}\S_{2})\circ (\tr\otimes \id)\circ \tr=0$,
\item $\im\tr\subset\ker (\Y_{2}(-\S))$.  \end{itemize} \end{DEF}

Now we consider the necessary and sufficient condition for the existence a
factor derivation of a factor algebra $\TG/\ker\Y (-\S)$.
\begin{LEM}[Identities for Y(-\S)]\label{iden} For the Young antisymmetrizer
$\Y_{3}(-\S)$ we have the following identity $$
\Y_{3}(-\S)\circ\S_{1}\circ\S_{2}=\Y_{3}(-\S)\circ\S_{2}\circ\S_{1}
=\Y_{3}(-\S).
$$
\end{LEM}

For the Young antisymmetrizer $\Y (-\S)$ a derivation of free tensor algebra
$\de\in\der(\TG)$ should satisfy condition imposing that it will be a factor
derivation $\deh\in\der (\TG/\ker \Y (-\S))$.  \begin{LEM}[Factor
derivation]\label{sfdlie} Let a comultiplication $\tr$ be $\S$-morphism and
for a derivation of a free algebra let $\de_{\tr}|\G=\tr$. Then exist a
factor derivation $\deh_{\tr}\in\der(\TG/\ker\Y(-\S))$, $$
(\tr\otimes\id)\S=\S_{1,2}(\id\otimes\tr) \qquad\Rightarrow\qquad
\de_{\tr}\ker\Y(-\S)\subset \ker\Y(-\S).  $$ \end{LEM} {\em Proof.} We can
write the condition $\de_{\tr}\ker\Y (-\S)\subset\ker\Y (-\S)$ in the
following form for a derivation $\de|\G^{\otimes 2}$ $$ \Y_{3} (-\S)\circ
 (\tr\otimes\id-\id\otimes\tr)\circ (\id+\S)=0.  $$ We can proof this
equation using the morphism condition and the lemma \ref{iden}. Using these
results we get $$[\Y_{3} (-\S)-\Y_{3}
(-\S)\circ\S_{1}\S_{2}]\circ(\tr\otimes\id)- [\Y_{3} (-\S)\circ\S_{2}\S_{1}-
\Y_{3} (-\S)]\circ(\id\otimes\tr)=0$$ $\eop$

The inverse is not true.  For the model (\ref{color}) of the symmetry $\S$
the condition $\de_{\tr}\ker\Y(-\S)\subset \ker\Y(-\S)$ are not equivalent to
the morphism condition for $\tr$.

\begin{LEM}[Identity for $\tr$ and $\S$]\label{blumen} Let the multiplication
$\tr$ be $\S$-morphism. Let $\tr_{\S}=\tr-\S\circ\tr$. Then we have the
following identity \begin{multline*} \Y (-\S)\circ
(\tr\otimes\id-\id\otimes\tr)\\
=[(\tr_{\S}\otimes\id)\circ\tr_{\S}-(\id\otimes\tr_{\S})\circ\tr_{\S}
-(\id\otimes\S)\circ(\tr_{\S}\otimes\id)\circ\tr_{\S}]\circ (\id-\S).
\end{multline*} \end{LEM} {\em Proof.} Using the $\S$-morphism condition for
the comultiplication and the ansatz with $\tr_{\S}$ above condition will be
obtained by simple calculations.  $\eop$

\begin{TH}[Koszul Theorem for $\S$-Lie coalgebras] For an $\S$-exterior
algebra ${\cal E}_{\S} (\G)$ over an S-Lie coalgebra $\SLA$ the following
assertions are equivalent \begin{itemize} \item $\deh^{2}=0\in
\G_{\S}^{\wedge 3},$ \item $\im (\de_{\tr})^{2}\subset\ker Y_{3} (-\S)$,
\item $(\id+\S_{1,2}+\S_{2,1})\circ(\id\otimes\tr)\circ\tr=0$, \item
$((\tr\otimes\id)\tr - (\id\otimes\tr)\circ\tr -
(\id\otimes\S)\circ(\tr\otimes\id)\circ\tr=0.$ \end{itemize} \end{TH} {\em
Proof.} If a comultiplication $\tr$ is $\S$-morphism then the first and the
second condition are equivalent from the lemma \ref{sfdlie}
$$\Y_{3} (-\S)\circ \de_{\tr}^{2}|\G=0\Leftrightarrow
\deh_{\tr}^{2}=0.$$ From the identity in lemma \ref{blumen} the second
condition is equivalent to the fourth condition, the Jacobi condition in the
Woronowicz form.  From this condition by the $\S$-commutativity of the
comultiplication $\tr$ we can get the third condition.  $\eop$

\section*{Acknowledgments} I would like to thank dr A. Borowiec, mgr M.
Ko\'{s}cielecki, prof. Z. Oziewicz and prof. S. Zakrzewski for remarks.  I am
also indebted to author wishes to thank dr B. Lulek and prof. T. Lulek for
the invitation and hospitality in Pozna\'{n}.


\begin{thebibliography}{10}

\bibitem{Boz} M. Bozejko and R. Speicher, \newblock Completely positive maps
on Coxeter groups, deformed commutation relations, and operator spaces,
\newblock Math. Ann. 300 (1994) 97--120.

\bibitem{DG} D. Gurevich, \newblock The Yang-Baxter equations and a
generlization of formal Lie theory, \newblock Soviet Math. Doklady 33 (1986)
758--762.

\bibitem{Kosz} J. Koszul, \newblock Homologie et cohomologie des alg\.ebras
de Lie, \newblock Bull. Soc. Math. France 78 (1950), 1--63.

\bibitem{SM} S. MacLane, \newblock Categories for the Working Mathematician,
Springer Verlag, 1971.

\bibitem{MR} N. Metropolis and G.-C. Rota, \newblock Symmetry Classes:
Function of Three Variables, \newblock American Mathematical Monthly 98 No. 4
(1991) 328--332.

\bibitem{WM} W. Michaelis, \newblock Lie coalgebras, \newblock Advances in
Mathematics 38 (1980), 1--54.

\bibitem{OPR} Z. Oziewicz, E. Paal and J. R\'{o}\.za\'{n}ski, \newblock
Coalgebras, cocompositions and cohomology, \newblock in `Non Associative
Algebra and Its Applications', \newblock in series Mathematics and its
Applications, ed. S. Gonzalez, Klu\-wer Acad. Publishers, str. 314--322,
1994.

\bibitem{Sch} Scheunert M., \newblock Generalized Lie algebras, \newblock J.
Math. Phys. 20 (1979).

\bibitem{SW} S. L. Woronowicz, \newblock Differential calculus on compact
matrix pseudogroups (quantum groups), \newblock Commun. Math. Phys. 122
(1989), 125--170.

\bibitem{DY} D. Yetter, \newblock Quantum groups and representations of
monoidal categories, \newblock Math. Proc. Camb. Phil. Soc. 108 (1990).

\end{thebibliography}
\end{document}